\documentclass[11pt]{amsart}

\usepackage{amsmath}
\usepackage{amssymb}
\usepackage{mathrsfs}
\usepackage{comment}
\usepackage{color}
\usepackage[colorlinks,citecolor=blue,urlcolor=black,linkcolor=black]{hyperref}

\newtheorem{theorem}{Theorem}[section]
\newtheorem{claim}[theorem]{Claim}
\newtheorem{mclaim}[theorem]{Main Claim}

\newtheorem{proposition}[theorem]{Proposition}

\theoremstyle{definition}
\newtheorem{definition}[theorem]{Definition}

\newtheorem{question}[theorem]{Question}

\theoremstyle{remark}
\newtheorem{remark}[theorem]{Remark}

\newcount\skewfactor
\def\mathunderaccent#1#2 {\let\theaccent#1\skewfactor#2
\mathpalette\putaccentunder}
\def\putaccentunder#1#2{\oalign{$#1#2$\crcr\hidewidth
\vbox to.2ex{\hbox{$#1\skew\skewfactor\theaccent{}$}\vss}\hidewidth}}
\def\name{\mathunderaccent\tilde-3 }

\def\smallbox#1{\leavevmode\thinspace\hbox{\vrule\vtop{\vbox
   {\hrule\kern1pt\hbox{\vphantom{\tt/}\thinspace{\tt#1}\thinspace}}
   \kern1pt\hrule}\vrule}\thinspace}

\newcommand{\cf}{{\rm cf}}

\def\qedref#1{$\qed_{\reforiginal{#1}}$}

\title{Magidor-like forcing and the cofinality of the Galvin number}

\author{Shimon Garti}
\address{Einstein Institute of Mathematics,
 The Hebrew University of Jerusalem,
 Jerusalem 91904, Israel}
\email{shimon.garty@mail.huji.ac.il}

\thanks{Research supported by ISF grant 2320/23.}

\begin{document}
\let\labeloriginal\label
\let\reforiginal\ref
\def\ref#1{\reforiginal{#1}}
\def\label#1{\labeloriginal{#1}}

\begin{abstract}
We prove that one can force over a model of \textsf{AD} to obtain countable cofinality of the full Galvin number. \\
Nous prouvons que on peut forcer sur un mod\`ele d'AD pour obtenir cofinalit\'e d\'enombrable du nombre de Galvin.
\end{abstract}

\maketitle

\newpage

\section{Introduction}

Let $\mathscr{F}$ be a normal filter over $\kappa$, where $\kappa$ is strongly regular and uncountable.\footnote{A cardinal $\kappa$ is \emph{strongly regular} if and only if $\kappa=\kappa^{<\kappa}$.}
Galvin proved that every family $\mathcal{C}=\{C_\alpha:\alpha\in\kappa^+\}\subseteq\mathscr{F}$ admits a subfamily $\mathcal{D}=\{C_{\alpha_i}:i\in\kappa\}$ so that $\bigcap\mathcal{D}\in\mathscr{F}$, see \cite{MR0369081}.
In particular, if $2^\omega=\omega_1$ then the club filter of $\aleph_1$ satisfies this property.

Given $\aleph_0<\kappa\leq\partial\leq\lambda$ and a normal filter $\mathscr{F}$ over $\kappa$ let us write ${\rm Gal}(\mathscr{F},\partial,\lambda)$ to denote the statement that every $\mathcal{C}\subseteq\mathscr{F}$ of size $\lambda$ admits a subfamily $\mathcal{D}\subseteq\mathcal{C}$ of size $\partial$ such that $\bigcap\mathcal{D}\in\mathscr{F}$.
Using this terminology, Galvin's theorem says that if $2^\omega=\omega_1$ then ${\rm Gal}(\mathscr{D}_{\aleph_1},\aleph_1,\aleph_2)$ holds, where $\mathscr{D}_{\aleph_1}$ is the club filter of $\aleph_1$.

Abraham and Shelah proved in \cite{MR830084} that one can force $2^{\omega_1}=\lambda$ along with $\neg{\rm Gal}(\mathscr{D}_{\aleph_1},\aleph_1,\lambda)$, where $\lambda$ is arbitrarily large.
An interesting feature of the proof is that $2^\omega=\lambda$ in this model as well.
It was shown in \cite{MR3604115} that this fact is not just a coincidence.
To wit, if $\mu=(2^\omega)^+$ then ${\rm Gal}(\mathscr{D}_{\aleph_1},\aleph_1,\mu)$ holds.
From the definition it follows that if $\mu'\geq\mu$ then ${\rm Gal}(\mathscr{F},\partial,\mu)$ implies ${\rm Gal}(\mathscr{F},\partial,\mu')$.
Thus, if $2^\omega<\lambda$ then ${\rm Gal}(\mathscr{D}_{\aleph_1},\aleph_1,\lambda)$ follows.
One concludes that Galvin's property for clubs of $\aleph_1$ depends on the size of $2^\omega$.
This fact is quite surprising, as one could expect that Galvin's property at $\aleph_1$ will be determined by the number of clubs of $\aleph_1$, that is $2^{\omega_1}$.
The fact that $2^\omega$ plays an important r\^ole here, gives rise to the following definition of $\mathfrak{gp}$ from \cite{MR3787522}, which places it in the area of cardinal characteristics of the continuum:

\begin{definition}
\label{defgp} The Galvin number. \newline
Let $\mathfrak{gp}$ be the minimal $\kappa$ such that every family $\mathcal{C}=\{C_\alpha:\alpha\in\kappa^+\}$ of clubs of $\aleph_1$ contains a subfamily $\mathcal{D}=\{C_{\alpha_i}:i\in\omega_1\}$ so that $\bigcap\mathcal{D}\in\mathscr{D}_{\aleph_1}$.
\end{definition}

A natural question is what the possible values of $\mathfrak{gp}$ are.
In the model of Abraham and Shelah, one begins with an infinite cardinal $\lambda$ and forces $2^\omega=2^{\omega_1}=\lambda$ along with the failure of ${\rm Gal}(\mathscr{D}_{\aleph_1},\aleph_1,\lambda)$.
Observe that ${\rm Gal}(\mathscr{D}_{\aleph_1},\aleph_1,(2^{\aleph_1})^+)$ is trivial, so ${\rm Gal}(\mathscr{D}_{\aleph_1},\aleph_1,\lambda^+)$ holds in this setting.
Consequently, $\mathfrak{gp}=\lambda$ in this model.
This forcing construction applies to every $\lambda$ so that $\cf(\lambda)>\aleph_1$, hence $\mathfrak{gp}$ can be either regular or singular of cofinality above $\aleph_1$, at every $\lambda$ in the $\aleph$-scale.

This pattern does not immediately work when $\cf(\lambda)\leq\aleph_1$, simply because $2^{\omega_1}=\lambda$ in the Abraham-Shelah model.
However, if $\cf(\lambda)=\omega_1$ then one can force $\mathfrak{gp}=\lambda$, essentially by the same forcing construction, as will be shown in the last section of this paper.
The case of countable cofinality is more interesting.
In particular, it was asked in \cite[Question 4.4]{MR3787522} whether $\cf(\mathfrak{gp})=\omega$ is consistent.
Though the question is still open, we know from \cite{MR4423479} that if $\cf(\mathfrak{gp})=\omega$ is consistent then it has some consistency strength.

The purpose of this paper is to build a model of \textsf{ZF} in which the cofinality of \emph{the full Galvin number} is countable.
The full Galvin number $\mathfrak{fgp}$ is a variation of $\mathfrak{gp}$, adapted to the context of \textsf{AD}.
A formal definition will be given in the next section.

The rest of the paper is arranged in three sections.
In the first section we give some background and in the second section we prove the main result.
In the last section we discuss some \textsf{ZFC} results with regard to the cofinality of the Galvin number.

Our notation is hopefully standard.
We shall employ the Jerusalem forcing notation, that is, we force upwards.
We denote the club filter over $\kappa$ by $\mathscr{D}_\kappa$ whenever $\kappa$ is regular and uncountable.
The collection of $\omega$-closed unbounded subsets of $\kappa$ is denoted by $\mathscr{U}^\kappa_\omega$.
In the case of $\kappa=\aleph_1$ this is simply the club filter $\mathscr{D}_{\aleph_1}$.
We let $\Theta$ be the supremum of $\lambda$ for which there is a surjection $\pi:{}^\omega\omega\rightarrow\lambda$.
Finally, a cardinal $\kappa$ is a boldface \textsf{GCH} cardinal iff there is no injection $f:\kappa^+\rightarrow\mathcal{P}(\kappa)$.

For background in \textsf{pcf} theory we refer to \cite{MR2768693} and to Shelah's monograph \cite{MR1318912}.
For background in the combinatorial aspects of \textsf{AD} we suggest the excellent monograph of Kleinberg, \cite{MR0479903}.
Basic information about polarized relations can be found in \cite{MR3075383}.

\newpage

\section{Preliminaries}

Seemingly, Galvin's property is connected to the axiom of choice.
The assumption $\kappa=\kappa^{<\kappa}>\aleph_0$ within the statement of Galvin's theorem requires some choice, and the proof meets this axiom as well.
However, instances of Galvin's property are provable even without full choice.
An initial study of this phenomenon is carried out in \cite{bgp}.
Here we phrase and prove several additional statements.

Ahead of our first claim, recall that $\binom{\alpha}{\beta}\rightarrow\binom{\gamma}{\delta}$ denotes the statement that for every coloring $c:\alpha\times\beta\rightarrow\{0,1\}$ one can find $A\in[\alpha]^\gamma,B\in[\beta]^\delta$ so that $c\upharpoonright(A\times{B})$ is constant.
One can replace the ordinals in this relation by some structure, with the natural intended meaning.
Thus, if $\mathscr{W}\subseteq\mathcal{P}(\lambda)$ and $\mathscr{U}\subseteq\mathcal{P}(\kappa)$ then $\binom{\lambda}{\kappa}\rightarrow\binom{\mathscr{W}}{\mathscr{U}}$ means that the sets $A,B$ which form the monochromatic product belong to $\mathscr{W}$ and $\mathscr{U}$ respectively.

\begin{claim} \label{clmmeasurables}
  Suppose that $\aleph_0<\kappa<\lambda$ and both $\kappa$ and $\lambda$ are measurable.
  Then ${\rm Gal}(\mathscr{U},\lambda,\lambda)$ holds whenever $\mathscr{U}$ is a normal ultrafilter over $\kappa$.
\end{claim}

\par\noindent\emph{Proof}. \newline
Let $\mathscr{U}$ be a normal ultrafilter over $\kappa$, and let $\mathscr{W}$ be a normal ultrafilter over $\lambda$.
As a first step we prove the combinatorial relation $\binom{\lambda}{\kappa}\rightarrow\binom{\mathscr{W}}{\mathscr{U}}$.
Then we derive the statement of the claim from this relation.

Assume, therefore, that $c:\lambda\times\kappa\rightarrow{2}$ is given.
For every $\beta\in\kappa$ there is a unique $j_\beta\in\{0,1\}$ for which $S^{j_\beta}_\beta=\{\alpha\in\lambda:c(\alpha,\beta)=j_\beta\}$ belongs to $\mathscr{W}$.
It follows that for some fixed $j$ and a set $B\in\mathscr{U}$ one has $j_\beta=j$ whenever $\beta\in{B}$.
Define $A=\bigcap\{S^j_\beta:\beta\in{B}\}$ and notice that $A\in\mathscr{W}$ by $\lambda$-completeness.
Observe that $c''(A\times{B})=\{j\}$, hence we are done proving that $\binom{\lambda}{\kappa}\rightarrow\binom{\mathscr{W}}{\mathscr{U}}$.

Assume now that $\{C_\alpha:\alpha\in\lambda\}\subseteq\mathscr{U}$.
Define $d:\lambda\times\kappa\rightarrow{2}$ by letting $d(\alpha,\beta)=0$ iff $\beta\in{C_\alpha}$.
Let $A\in\mathscr{W},B\in\mathscr{U}$ be such that $d\upharpoonright(A\times{B})$ is constant.
Notice that the constant value cannot be $1$.
Indeed, choose $\alpha\in{A}$ and any element $\beta\in B\cap C_\alpha$ and conclude that $\beta\notin C_\alpha$ if $d''(A\times{B})=\{1\}$, contradicting the choice of $\beta$.
Thus $d''(A\times{B})=\{0\}$, and hence $B\subseteq\bigcap\{C_\alpha:\alpha\in{A}\}$.
Since $B\in\mathscr{U}$, the proof is accomplished.

\hfill \qedref{clmmeasurables}

Let us indicate that in order to obtain the statement ${\rm Gal}(\mathscr{U},\lambda,\lambda)$ from the relevant polarized relation one has to assume that $\kappa$, the domain of $\mathscr{U}$, is measurable, but the assumption that $\lambda$ is also measurable can be relaxed.
A moment of perusal shows that if $\binom{\lambda}{\kappa}\rightarrow\binom{\lambda}{\mathscr{U}}$ then ${\rm Gal}(\mathscr{U},\lambda,\lambda)$ obtains.
The argument is just the same, and this fact might be helpful.

The statement ${\rm Gal}(\mathscr{U},\lambda,\lambda)$ is a strong form of the Galvin property.
The above claim establishes a strong Galvin property which will be dubbed as \emph{full}.
In \textsf{ZFC} this claim is trivial (if $\kappa<\lambda$ and both are measurables) since $\lambda=\cf(\lambda)>2^\kappa$, but without full choice this claim is informative.
Under \textsf{AD} it is frequent to meet this strong form, since the distance between regular cardinals and measurable cardinals below $\Theta$ is very small.\footnote{Actually, if $V=L(\mathbb{R})$ then every regular cardinal below $\Theta$ is measurable under \textsf{AD}, see \cite{MR2768698}.}
We define, therefore, the following variant of the Galvin number:

\begin{definition} \label{deffgp}
  The full Galvin number. \newline
  Let $\mathfrak{fgp}$ be the first $\lambda$ for which every family $\mathcal{C}=\{C_\alpha:\alpha\in\lambda^+\}\subseteq\mathscr{D}_{\aleph_1}$ contains a subfamily $\mathcal{D}$ of size $\lambda^+$ such that $\bigcap\mathcal{D}\in\mathscr{D}_{\aleph_1}$.
\end{definition}

By Claim \ref{clmmeasurables} and the fact that both $\aleph_1$ and $\aleph_2$ are measurable cardinals under \textsf{AD} we see that $\textsf{AD}\vdash\mathfrak{fgp}=\aleph_1$.
Our main goal is to force $\mathfrak{fgp}=\aleph_\omega$ in some model of \textsf{ZF}.
But let us prove, first, that the statement $\cf(\mathfrak{fgp})=\omega$ has some consistency strength in \textsf{ZFC}.
The argument is similar to the parallel argument with respect to $\mathfrak{gp}$ as given in \cite{MR4423479}, but the details are somewhat different.

\begin{proposition}[ZFC] \label{propfgpstrength}
  The statement $\cf(\mathfrak{fgp})=\omega$ has some consistency strength.
\end{proposition}

\par\noindent\emph{Proof}. \newline
The proof consists of two steps.
In the first step we show that if $\mu$ is a singular cardinal of countable cofinality, $\mu=\bigcup_{n\in\omega}\mu_n$ where $(\mu_n:n\in\omega)$ is an increasing sequence of regular cardinals, ${\rm tcf}(\prod_{n\in\omega}\mu_n,J^{\rm bd}_\omega)=\mu^+$ and $(f_\alpha:\alpha\in\mu^+)$ is a scale witnessing this fact, endowed with the property that for every $A\in[\mu^+]^{\aleph_1}$ one has $|\{f_\alpha(n):\alpha\in{A}, n\in\omega\}|=\aleph_1$ then $\mathfrak{fgp}\neq\mu$.

By way of contradiction assume that the above assumptions hold yet $\mathfrak{fgp}=\mu$.
Thus, for every $n\in\omega$ one has $\mathfrak{fgp}>\mu_n$ and hence one can choose a family $\mathcal{D}_n=\{D^n_\alpha:\alpha\in\mu_n\}\subseteq\mathscr{D}_{\aleph_1}$ such that $\bigcap\{D^n_{\alpha_i}:i\in{B}\}$ is bounded in $\omega_1$ whenever $B\in[\mu_n]^{\aleph_1}$.
Set $\mathcal{D}=\bigcup\{\mathcal{D}_n:n\in\omega\}$ and notice that $|\mathcal{D}|=\mu$.
If $B'\in[\mu]^{\aleph_1}$ then there must be a set $B\in[B']^{\aleph_1}$ and a natural number $n$ such that $B\in[\mu_n]^{\aleph_1}$.
One concludes now that:
\[|\bigcap\{D^n_{\alpha}:\alpha\in{B'},n\in\omega\}|\leq|\bigcap\{D^n_{\alpha_i}:i\in{B}\}|<\aleph_1\]
This is the property that we need with respect to the family $\mathcal{D}$.

For every $\alpha\in\mu^+$ let $C_\alpha=\bigcap\{D^n_{f_\alpha(n)}:n\in\omega\}$.
Observe that $C_\alpha\in\mathscr{D}_{\aleph_1}$, being an intersection of but $\aleph_0$ elements of $\mathscr{D}_{\aleph_1}$.
Let $\mathcal{C}=\{C_\alpha:\alpha\in\mu^+\}$.
One can verify that $|\mathcal{C}|=\mu^+$ since $(f_\alpha:\alpha\in\mu^+)$ is a scale.
By our assumption, $\mathfrak{fgp}=\mu$.
Therefore, one can find $A\in[\mu]^{\aleph_1}$ such that the set $C=\bigcap\{C_\alpha:\alpha\in{A}\}$ belongs to $\mathscr{D}_{\aleph_1}$.
By our assumptions, $|\{f_\alpha(n):\alpha\in{A}, n\in\omega\}|=\aleph_1$.
Let $B=\{f_\alpha(n):\alpha\in{A},n\in\omega\}$, so $B\in[\mu]^{\aleph_1}$.
In particular, the set $D=\bigcap\{D^n_{f_\alpha(n)}:\alpha\in{A},n\in\omega\}$ is bounded in $\omega_1$.
However, $C=D$ by the above definitions, a contradiction.

We proceed now to the second step of the proof.
Suppose that $\mu>\cf(\mu)=\omega$ and there is an increasing sequence $(\mu_n:n\in\omega)$ of regular cardinals such that $\mu=\bigcup_{n\in\omega}\mu_n$ and ${\rm tcf}(\prod_{n\in\omega}\mu_n,J^{\rm bd}_\omega)=\mu^+$.
Assume, further, that $(f_\alpha:\alpha\in\mu^+)$ is a \emph{good scale} witnessing this fact.
From \cite[Theorem 1.4]{MR4423479} we know that if $A\in[\mu^+]^{\aleph_1}$ then $|\{f_\alpha(n):\alpha\in{A},n\in\omega\}|=\aleph_1$.
It follows now from the first step of the proof that $\mathfrak{fgp}\neq\mu$.

Assume that $\square^*_\mu$ holds at every $\mu\in{2^\omega}$ such that $\mu>\cf(\mu)=\omega$.
It is known that under this assumption, for every such $\mu$ there exists an increasing sequence of regular cardinals $(\mu_n:n\in\omega)$ such that $\mu=\bigcup_{n\in\omega}\mu_n$ and ${\rm tcf}(\prod_{n\in\omega}\mu_n,J^{\rm bd}_\omega)=\mu^+$ and there is a good scale witnessing the latter fact.\footnote{See, e.g. \cite[Section 15]{MR2160657}.}
Consequently, $\mathfrak{fgp}\neq\mu$ for every such $\mu$ and hence, in particular, $\cf(\mathfrak{fgp})>\omega$.
It follows that if one wishes to force $\cf(\mathfrak{fgp})=\omega$ then one must force $\neg\square^*_\mu$ for every $\mu$ of countable cofinality below the continuum, and this setting has some consistency strength, so we are done.

\hfill \qedref{propfgpstrength}

Back to models of \textsf{ZF}, our plan is to begin with a model of \textsf{AD} and to force over it.
The generic extension will not be a model of \textsf{AD} anymore, but it will be a model of \textsf{ZF} in which $\mathfrak{fgp}=\aleph_\omega$.
We shall force with a Magidor-like forcing to singularize $\aleph_2$, and the cofinality target is $\aleph_1$.
The details of this forcing notion are defined and proved in \cite{MR722169}.

Let us indicate that this forcing notion is similar to Prikry forcing, and differs from the classical (and more complicated) Magidor forcing as articulated in \cite{MR465868} in the context of \textsf{ZFC}.
The main challenge with forcing a measurable cardinal to be singular with uncountable cofinality is to prove the Prikry property.
In the case of Prikry forcing and countable cofinality, this property follows from a powerful partition theorem which holds at measurable cardinals.
But in the case of Magidor forcing and uncountable cofinality, the pertinent statement should be a partition theorem of infinite subsets of the measurable cardinal, and this is not available in \textsf{ZFC}.
Thus, Magidor forcing is an intricate combination of Prikry sequences which are forced into measurable cardinals below the measurable cardinal to be singularized with uncountable cofinality.
In particular, one adds many bounded subsets to $\kappa$ while performing this process.

Working in \textsf{ZF}, one can make strong assumptions about the measurable cardinal to be singularized, and partition theorems for infinite subsets of this cardinal are at hand.
This simplifies the definition of Magidor forcing considerably.
For example, in order to force cofinality $\omega_1$ one can force with countable conditions, much like Prikry forcing.
In particular, Magidor-like forcing in this formulation adds no bounded subsets of the given measurable cardinal.

We shall use the notation of Henle in \cite{MR722169}, and we mention some conventions.
Let $\kappa$ be a regular and uncountable cardinal.
For $A\in[\kappa]^\theta$ let $A_\alpha$ be the $\alpha$th element of $A$ and let ${}_\omega{A}$ be the increasing enumeration of countable sups of elements of $A$, that is $\{\bigcup_{n\in\omega}A_{\alpha+n}:\alpha\in\theta\}$. If $A\in[\kappa]^\kappa$ then let $\langle{A}\rangle$ be the set $\{{}_\omega B:B\in[A]^\kappa\}$.
We shall force with the following:

\begin{definition} \label{defmagidorlike} The forcing $\mathbb{M}_{\gamma\kappa}$. \newline
  Let $\kappa$ be regular and uncountable and assume that $\omega\leq\gamma<\kappa$.
  \begin{enumerate}
    \item [$(\aleph)$] A condition $p\in\mathbb{M}_{\gamma\kappa}$ is a pair $(s^p,A^p)$ such that $s^p\in[\kappa]^{<\gamma}, A^p\in[\kappa]^\kappa$ and $\bigcup{s^p}<\bigcap{A^p}$.
    \item [$(\beth)$] If $p,q\in\mathbb{M}_{\gamma\kappa}$ then $p\leq{q}$ iff $s^p\subseteq{s^q}, \langle{A^q}\rangle\subseteq\langle{A^p}\rangle$ and $s^q-s^p={}_\omega{t}$ for some $t\in[A]^{<\gamma}$.
  \end{enumerate}
\end{definition}

The definition is given in a general form, but some assumptions about $\kappa$ must be made so that the forcing will not collapse cardinals.
We shall say that $\kappa$ is a \emph{weak partition cardinal} if $\kappa\rightarrow(\kappa)^\gamma$ for every $\gamma\in\kappa$.
Henle proved in \cite{MR722169} that if $\kappa$ is a weak partition cardinal and $\omega<\gamma<\kappa$ then $\mathbb{M}_{\gamma\kappa}$ is a forcing notion which satisfies the Prikry property and preserves cardinals.
Moreover, if $G\subseteq\mathbb{M}_{\gamma\kappa}$ is generic then $V[G]$ contains no new bounded subsets of $\kappa$.

\begin{claim} \label{clmpreservation}
  Let $V$ be a model of $\mathsf{AD}$ and let $G\subseteq\mathbb{M}_{\aleph_1\aleph_2}$ be $V$-generic.
  Then $\aleph_1$ is measurable in $V[G]$ and $\mathscr{D}_{\aleph_1}$ is an ultrafilter.
\end{claim}

\par\noindent\emph{Proof}. \newline
By a celebrated theorem of Solovay, $\aleph_1$ is measurable in $V$ as witnessed by $\mathscr{D}_{\aleph_1}$.
Kleinberg proved in \cite{MR0479903} that $\aleph_2$ is a weak partition cardinal under \textsf{AD}.
Thus, $\mathbb{M}_{\aleph_1\aleph_2}$ satisfies the above mentioned properties and, in particular, adds no subsets of $\aleph_1$ since any such set is bounded in $\aleph_2$.
It follows that $\mathscr{D}_{\aleph_1}$ is still a normal ultrafilter over $\aleph_1$ in the generic extension by $\mathbb{M}_{\aleph_1\aleph_2}$, so we are done.

\hfill \qedref{clmpreservation}

Let us indicate that in the definition of the forcing order, $p\leq{q}$ contains the proviso that $s^p\subseteq{s^q}$.
We emphasize that the relation here is \emph{inclusion} and not \emph{end-extension} as appears usually in Prikry-type forcing notions.
This point plays an important role in some arguments below.

\newpage

\section{Countable cofinality}

``But it was pretty and characteristic, besides being singular" (\cite[Chapter IV]{tale}).
In this section we prove the central result of the paper.
We need a few facts ahead of the proof.
The first one, dictated in the main claim below, is a specific case of a result of Henle from \cite{MR722169}.
For completeness, we unfold the proof.

\begin{mclaim} \label{mclmpreservation}
  Assume $\mathsf{AD}$ and let $G\subseteq\mathbb{M}_{\aleph_1\aleph_2}$ be $V$-generic.
  Then $\aleph_{\omega+1}$ is regular in $V[G]$.
\end{mclaim}

\par\noindent\emph{Proof}. \newline
Recall that $\cf(\aleph_{\omega+1})^V=\aleph_{\omega+1}>\aleph_2$, since $\aleph_{\omega+1}$ is a measurable cardinal in $V$.
Fix a $V$-generic set $G\subseteq\mathbb{M}_{\aleph_1\aleph_2}$.
Our goal is to show that $\cf(\aleph_{\omega+1})^{V[G]}=\aleph_{\omega+1}$ as well.

Assume towards contradiction that there is a cardinal $\lambda\in\aleph_{\omega+1}$ and a condition $(a,p)\in\mathbb{M}_{\aleph_1\aleph_2}$ such that $(a,p)\Vdash f:\lambda\rightarrow\aleph_{\omega+1}$, $f$ is order preserving and $\bigcup f''\lambda=\aleph_{\omega+1}$.
For every $s\in[p]^{<\omega_1}$ and every $\alpha\in\lambda$ we define the ordinal $\eta_{s\alpha}$ as follows.
If there exists $t\in[\omega_2]^{\omega_2}$ such that $(a\cup{}_\omega s,t)\Vdash\name{f}(\alpha)=\check{\beta}$ then $\eta_{s\alpha}=\beta$.
If not, let $\eta_{s\alpha}=0$.
Notice that $\eta_{s\alpha}$ is well-defined.
Indeed, if $t_0,t_1\in[\omega_2]^{\omega_2}$ and $(a\cup{}_\omega s,t_0)\Vdash\name{f}(\alpha)=\check{\beta_0},(a\cup{}_\omega s,t_1)\Vdash\name{f}(\alpha)=\check{\beta_1}$ then necessarily $\beta_0=\beta_1$ since conditions with the same stem are compatible.
Now for every $s\in[p]^{<\omega_1}$ let $H_s=\{\eta_{s\alpha}:\alpha\in\lambda\}$, and let $\eta_s=\bigcup H_s$, thus $\eta_s\leq\lambda$.

Define a coloring $c:[p]^{<\omega_1}\rightarrow{}^{\omega_1}2$ as follows.
Given $s\in[p]^{<\omega_1}$ and $\sigma\in\omega_1$, let $c(s)_\sigma=1$ iff there are $\sigma_0,\sigma_1\in\omega_1$ so that $\sigma<\sigma_0<\sigma_1$ and $\eta_{s\upharpoonright\sigma}<\eta_{s\upharpoonright(\sigma_1-\sigma_0)}$.
Recall that $\aleph_2\rightarrow(\aleph_2)^{<\omega_1}$ in $V$, and let $q\in[p]^{\omega_2}$ be $c$-monochromatic.

The crucial point here is that for every $\sigma\in\omega_1$ one can find $\sigma'>\sigma, \sigma'\in\omega_1$ and $s'\in[q]^{\sigma'}$ for which $c(s')_{\sigma'}=1$.
To see this, fix $\sigma\in\omega_1$ and let $s=q\upharpoonright\sigma$.
Pick $r\in[q-s]^{<\omega_1}, t\in[q-s]^{<\omega_2}$ and $\alpha\in\lambda$ such that $(a\cup{}_\omega r,t)\Vdash\name{f}(\alpha)=\beta>\eta_s$.
This is possible since $q\subseteq p$ and $(a,p)\Vdash\bigcup f''\lambda=\aleph_{\omega+1}$.
Observe that $c(s\cup r)_\sigma=1$ since $\eta_r\geq\beta>\eta_s$.
Let $s'=s\cup r$ and let $\sigma'=\ell g(s')$, thus the crucial property holds.

Suppose that $s$ is a sequence of \emph{consecutive elements} from $q$ of length $\tau$.
We shall say that $s$ is $\tau$-good (or $\tau$-good for $q$), and denote it by $s^g_\tau$.
Define $R=\{\eta_{s^g_\tau\alpha}:\alpha\in\lambda,\tau\in\omega_1\}$.
Notice that $|R|<\aleph_{\omega+1}$.
However, $(a,q)\Vdash f''\lambda\subseteq{R}$, since $q$ is $c$-monochromatic, and hence $f''\lambda$ is bounded in $\aleph_{\omega+1}$.
This is a contradiction, so the proof is accomplished.

\hfill \qedref{mclmpreservation}

We saw already that positive polarized relations are deeply connected with Galvin's property.
This can be seen in both directions.
We shall use the following negative relation in order to show that the Galvin property fails in the relevant context.

\begin{claim}
  \label{clmsamecof}
  \begin{enumerate}
    \item [$(a)$] If $\kappa=\cf(\lambda)$ then $\binom{\lambda}{\kappa}\nrightarrow\binom{\lambda}{\kappa}$.
    \item [$(b)$] If ${\rm Gal}(\mathscr{D}_\kappa,\lambda,\lambda)$ holds then $\binom{\lambda}{\kappa}\rightarrow\binom{\lambda}{\kappa}$.
    \item [$(c)$] If $\kappa=\cf(\lambda)$ then $\neg{\rm Gal}(\mathscr{D}_\kappa,\lambda,\lambda)$.
  \end{enumerate}
\end{claim}

\par\noindent\emph{Proof}. \newline
We commence with part $(a)$.
Fix an increasing and continuous sequence $(\lambda_i:i\in\kappa)$ of ordinals such that $\lambda=\bigcup_{i\in\kappa}\lambda_i$.
For every $\alpha\in\lambda$ let $\eta(\alpha)$ be the unique ordinal $i\in\kappa$ so that $\alpha\in[\lambda_i,\lambda_{i+1})$.
Define a coloring $c:\lambda\times\kappa\rightarrow{2}$ as follows.
Given $\alpha\in\lambda,\beta\in\kappa$ let $c(\alpha,\beta)=0$ iff $\eta(\alpha)\leq\beta$.

Suppose that $A\in[\lambda]^\lambda,B\in[\kappa]^\kappa$.
Fix $\alpha\in{A}$.
Choose $\beta\in{B}$ so that $\eta(\alpha)\leq\beta$.
There must be such an ordinal since $B$ is unbounded in $\kappa$.
By definition, $c(\alpha,\beta)=0$.
Now pick $\alpha'\in{A}$ so that $\eta(\alpha')>\beta$.
There must be such an ordinal since $A$ is unbounded in $\lambda$.
By definition, $c(\alpha',\beta)=1$.
We conclude that $c''(A\times{B})=\{0,1\}$, as required.

For part $(b)$, assume that ${\rm Gal}(\mathscr{D}_\kappa,\lambda,\lambda)$ holds and let $c:\lambda\times\kappa\rightarrow{2}$ be any coloring.
For every $\alpha\in\lambda$ there is (a unique) $i_\alpha\in\{0,1\}$ such that $S^{i_\alpha}_\alpha=\{\beta\in\kappa:c(\alpha,\beta)=i_\alpha\}\in\mathscr{D}_\kappa$.
Fix $A'\in[\lambda]^\lambda,i\in\{0,1\}$ so that $i_\alpha=i$ whenever $\alpha\in{A'}$.
Let $\mathcal{D}=\{S^i_\alpha:\alpha\in{A'}$.
From ${\rm Gal}(\mathscr{D}_\kappa,\lambda,\lambda)$ we infer that there are $B\in\mathscr{D}_\kappa$ and $A\in[A']^\lambda$ such that $B\subseteq\bigcap\{S^i_\alpha:\alpha\in{A}$.
By definition, $c''(A\times{B})=\{i\}$, so we are done proving $(b)$.
Finally, part $(c)$ follows from parts $(a)$ and $(b)$.

\hfill \qedref{clmsamecof}

Equipped with the above claims we can prove the following:

\begin{theorem} \label{thmcffgp}
  One can force over a model of $\mathsf{AD}$ to obtain $\mathfrak{fgp}=\aleph_\omega$.
\end{theorem}

\par\noindent\emph{Proof}. \newline
Let $V$ be a model of \textsf{AD}, and assume further that $V=L(\mathbb{R})$.
Recall that $\aleph_1, \aleph_2$ and $\aleph_{\omega+1}$ are measurable in $V$, and in particular $\aleph_{\omega+1}$ is regular.
Let $G\subseteq\mathbb{M}_{\aleph_1\aleph_2}$ be $V$-generic, so $V[G]\models\cf(\aleph_2)=\omega_1$, thus $V[G]$ is not a model of \textsf{AD}.
However, $\mathscr{D}_{\aleph_1}$ is still a normal ultrafilter in $V[G]$ by virtue of Claim \ref{clmpreservation}.
Likewise, $V[G]\models\aleph_{\omega+1}$ is a regular cardinal due to the main claim.

Suppose that $\langle A_\alpha:\alpha\in\aleph_{\omega+1}\rangle$ is a sequence of elements of $\mathscr{D}_{\aleph_1}$ in $V[G]$.
Since $V[G]$ contains no new subsets of $\aleph_1$ one concludes from \cite[Theorem 8.26]{MR2768698} that $|\{A_\alpha:\alpha\in\aleph_{\omega+1}\}|\leq\aleph_1$.
Since $\aleph_{\omega+1}$ is regular in $V[G]$, there are $T\in[\aleph_{\omega+1}]^{\aleph_{\omega+1}}$ and $A\in\mathscr{D}_{\aleph_1}$ such that $\alpha\in{T}\Rightarrow A_\alpha=A$.
This fact witnesses ${\rm Gal}(\mathscr{D}_{\aleph_1},\aleph_{\omega+1},\aleph_{\omega+1})$ and hence $\mathfrak{fgp}\leq\aleph_\omega$ in $V[G]$.

Trivially, $\mathfrak{fgp}\geq\aleph_1$, e.g. by the collection of end-segments of $\aleph_1$.
So let us prove that if $0<n<\omega$ then $\mathfrak{fgp}\neq\aleph_n$.
In other words, let us show that $\neg{\rm Gal}(\mathscr{D}_{\aleph_1},\aleph_{n+1},\aleph_{n+1})$ whenever $n<\omega$.

Recall that $\cf(\aleph_n)=\omega_2$ under \textsf{AD} for every $n\in[2,\omega)$.
Hence in $V[G]$ one has $\cf(\aleph_n)=\omega_1$ for every $n\in[2,\omega)$.
Similarly, $\cf(\aleph_1)=\aleph_1$ in $V[G]$, since $\aleph_1$ is measurable and hence regular.
Consequently, if $n>1$ then $\binom{\aleph_n}{\aleph_1}\nrightarrow\binom{\aleph_n}{\aleph_1}$ by virtue of Claim \ref{clmsamecof}(a).
From Claim \ref{clmsamecof}(c) we infer that if $0<n<\omega$ (so $1<n+1<\omega$) then $\neg{\rm Gal}(\mathscr{D}_{\aleph_1},\aleph_{n+1},\aleph_{n+1})$ and hence $\mathfrak{fgp}\neq\aleph_n$.
Therefore, $\mathfrak{fgp}\geq\aleph_\omega$ and then $\mathfrak{fgp}=\aleph_\omega$ as required.

\hfill \qedref{thmcffgp}

It is possible, of course, to force Prikry into $\aleph_2$, and then the cofinality of every $\aleph_n$ above $\aleph_1$ will be countable.
However, we do not know whether $\mathfrak{fgp}=\aleph_\omega$ in this model.
Polarized relations of the form $\binom{\mu}{\kappa}\rightarrow\binom{\mu}{\kappa}$ where $\kappa=\cf(\kappa)>\mu$ and $\mu>\cf(\mu)$ were studied in \cite{MR0202613}.
The authors were interested, inter alia, in the situation in which $\kappa=\cf(\mu)^+$.
It is shown in \cite{MR4101445} that such relations are independent.

In our context, we need the negative relation $\binom{\aleph_n}{\aleph_1}\nrightarrow\binom{\aleph_n}{\aleph_1}$ for every $n\in[2,\omega)$.
After Prikry forcing, the cofinality of $\aleph_n$ for $n\in[2,\omega)$ becomes $\aleph_0$, and it is not clear whether the above negative relation holds.
Notice that $\binom{\aleph_1}{\aleph_0}\rightarrow\binom{\aleph_1}{\aleph_0}$ under \textsf{AD}, so $\binom{\aleph_n}{\aleph_1}\rightarrow\binom{\aleph_n}{\aleph_1}$ is quite reasonable in this context.
For that reason, we force with Magidor forcing into $\aleph_2$ to get $\cf(\aleph_2)=\omega_1$ in the generic extension.

In the proof of the above theorem we used the fact that $\aleph_{\omega+1}$ is still regular in the generic extension by $\mathbb{M}_{\aleph_1\aleph_2}$.
There is another way to achieve the same result, provided that $\aleph_{\omega+1}$ remains measurable.
The advantage of this approach is that the (seemingly stronger) assumption of measurability yields the pertinent Galvin property even if $\aleph_1$ is not boldface \textsf{GCH}.
Thus, if one forces over a model of \textsf{AD} (in which $\aleph_{\omega+1}$ is measurable) then one should ask whether $\aleph_{\omega+1}$ remains measurable in the generic extension by $\mathbb{M}_{\aleph_1\aleph_2}$.

A celebrated theorem of L\'evy and Solovay from \cite{MR224458} says that if $\kappa$ is measurable and $\mathbb{P}$ is a forcing notion of size less than $\kappa$ then $\kappa$ remains measurable in any generic extension by $\mathbb{P}$.
We need a similar preservation theorem, but we work in the context of \textsf{AD} in which in most cases (including the case of Magidor-like forcing) we cannot apply an assumption like $|\mathbb{P}|<\kappa$.

Nonetheless, it seems reasonable to assume that the measurability of $\aleph_{\omega+1}$ will be preserved if one forces with $\mathbb{M}_{\aleph_1\aleph_2}$ over a model of \textsf{AD}.
We indicate that if one forces with Prikry forcing into $\aleph_1$ over a model of \textsf{AD} then $\aleph_2$ remains measurable in the generic extension,\footnote{This is a result of Henle, see \cite[Proposition 3.5]{MR722169}.} so the parallel situation with respect to $\aleph_{\omega_1}$ is plausible.
There is a meaningful difference, however, between these two cases.
In the case of $\aleph_1$ and $\aleph_2$, the latter is isomorphic to ${}^{\omega_1}\omega_1/\mathscr{U}^{\omega_1}_\omega$, and this fact plays an important role in the proof that the measurability of $\aleph_2$ is preserved, see \cite[Proposition 3.5]{MR722169}.
In the case of $\aleph_2$ and $\aleph_{\omega+1}$, this is not the case anymore.

If, however, $\aleph_{\omega+1}$ remains measurable, then $\binom{\aleph_{\omega+1}}{\aleph_1}\rightarrow\binom{\mathscr{W}}{\mathscr{D}_{\aleph_1}}$, where $\mathscr{W}$ is a normal ultrafilter over $\aleph_{\omega+1}$, as proved in Claim \ref{clmmeasurables}.
From this one concludes that ${\rm Gal}(\mathscr{D}_{\aleph_1},\aleph_{\omega+1},\aleph_{\omega+1})$ holds and hence $\mathfrak{fgp}\leq\aleph_\omega$ in $V[G]$.
The opposite direction is proved exactly as in the proof of the above theorem.
We are left, therefore, with the following:

\begin{question}
  \label{qmeasurability} Let $V$ be a model of \textsf{AD}, and let $G$ be $V$-generic for $\mathbb{M}_{\aleph_1\aleph_2}$.
  \begin{enumerate}
    \item [$(\aleph)$] Is $\aleph_{\omega+1}$ a measurable cardinal in $V[G]$?
    \item [$(\beth)$] Is the collection of $\omega$-closed unbounded subsets of $\aleph_{\omega+1}$ a normal ultrafilter in $V[G]$?
  \end{enumerate}
\end{question}

In this section we forced $\mathfrak{gp}=\aleph_\omega$.
We started with a model of \textsf{AD} in which $\aleph_{\omega+1}$ is measurable, and utilized Magidor forcing $\mathbb{M}_{\aleph_1\aleph_2}$.
Let us observe that the value of $\mathfrak{gp}$ is not confined to $\aleph_\omega$ using this method.
That is, one can force $\cf(\mathfrak{gp})=\omega$ with larger values of $\mathfrak{gp}$.

For example, one can start with $\mathbb{M}_{\aleph_1\kappa}$ where $\kappa=\aleph_{\omega+2}$ followed by the same forcing notion where $\kappa=\aleph_{\omega+1}$.
Since Magidor forcing in this context adds no bounded subsets to $\kappa$, $\aleph_{\omega+1}$ remains measurable after the first step, and hence one can perform the second step.
Finally, one should force with $\mathbb{M}_{\aleph_1\aleph_2}$, using again the fact that $\aleph_2$ will remain measurable after the first two steps of the iteration.
In the resulting model, all the cardinals up to $\aleph_{\omega^{\omega^\omega}+1}$ are either of countable cofinality or of cofinality $\aleph_1$.
On the other hand, $\aleph_{\omega^{\omega^\omega}+1}$ remains a regular cardinal, by the arguments of Main Claim \ref{mclmpreservation}.
Thus, $\mathfrak{gp}=\aleph_{\omega^{\omega^\omega}+1}$ in the generic extension.

\newpage

\section{Remarks on the Galvin number in ZFC}

In this section we discuss two issues related to $\mathfrak{gp}$ in \textsf{ZFC}.
First, we observe that one of the statements of \cite{MR4423479} is incorrect, and we point out the mistake.
Moreover, we prove here the opposite statement.
Second, we indicate that if $\cf(\mathfrak{gp})=\omega$ then a specific constellation of cardinal arithmetic must hold.

Suppose that $\theta=\cf(\theta)>\aleph_1$.
One can force $\mathfrak{gp}=\theta$ over a model of \textsf{GCH}, thus $\cf(\mathfrak{gp})=\theta$ is consistent.
Moreover, if $\lambda>\cf(\lambda)=\theta$ then the model of \cite{MR830084} (with $\lambda$ as a parameter) gives $\mathfrak{gp}=\lambda$, so the values of $\cf(\mathfrak{gp})$ under question are just $\aleph_0$ and $\aleph_1$.

It was claimed in \cite{MR4423479} that both statements $\cf(\mathfrak{gp})=\aleph_0$ and $\cf(\mathfrak{gp})=\aleph_1$ have some consistency strength.
It seems, however, that the argument for $\cf(\mathfrak{gp})=\aleph_1$ is false.
Let us explain, briefly, the problematic point.
The explanation will be followed by a proof of the consistency of $\cf(\mathfrak{gp})=\aleph_1$, without any use of large cardinals.

Fix $\mu>\cf(\mu)=\omega$.
If $\mathfrak{gp}=\mu$ then there is a family $\mathcal{D}=\{D_\alpha:\alpha\in\mu\}\subseteq\mathscr{D}_{\aleph_1}$ such that $\bigcap\{D_\alpha:\alpha\in{B}\}$ is not in $\mathscr{D}_{\aleph_1}$ whenever $B\in[\mu]^{\aleph_1}$.
Indeed, choose an increasing sequence $(\mu_n:n\in\omega)$ of regular cardinals such that $\mu=\bigcup_{n\in\omega}\mu_n$.
For every $n\in\omega$ we are assuming that $\mathfrak{gp}>\mu_n$.
Hence there is a family $\mathcal{D}_n=\{D^n_\alpha:\alpha\in\mu_n^+\}$ witnessing this fact.
Let $\mathcal{D}=\bigcup_{n\in\omega}\mathcal{D}_n$.
Now if $B\in[\mu]^{\aleph_1}$ then there are $n\in\omega$ and $C\in[B]^{\aleph_1}$ so that $C\in[\mu_n^+]^{\aleph_1}$.
It follows that $D=\bigcap\{D_\alpha:\alpha\in{C}\}$ is not in $\mathscr{D}_{\aleph_1}$ by the choice of $\mathcal{D}_n$.

This argument breaks down if $\mu>\cf(\mu)=\omega_1$.
One can choose a sequence $(\mu_i:i\in\omega_1)$ and $\mathcal{D}_i$ for every $i\in\omega_1$ as before.
But now it is possible that $B\in[\mu]^{\aleph_1}$ while $|B\cap\mu_i|<\aleph_1$ for each $i\in\omega_1$.
Therefore, the argument for the consistency strength of the statement $\cf(\mathfrak{gp})=\omega_1$ is false.

In fact, we shall prove below that one can force $\mathfrak{gp}=\lambda$ where $\lambda>\cf(\lambda)=\omega_1$, and we do not need large cardinals in the ground model for this statement.
We shall use the following terminology from \cite{MR4611828}.
Suppose that $\mathcal{C}=\{C_\alpha:\alpha\in\lambda\}\subseteq\mathscr{D}_{\aleph_1}$ witnesses the failure of the Galvin property.
We shall say that $\mathcal{C}$ witnesses the \emph{ultimate failure of the Galvin property} if $\lambda=2^{\omega_1}$.
This means that the size of the witness is the largest possible size for such a family.
If $\cf(\lambda)=\omega_1$ then $2^{\omega_1}\neq\lambda$, so if $\mathfrak{gp}=\mu$ exemplifies the ultimate failure of Galvin's property then $\cf(\mu)>\omega_1$.
Nevertheless, $\mathfrak{gp}=\lambda<2^{\omega_1}$ is possible even if $\cf(\lambda)=\omega_1$, as shown in the following.

\begin{claim}
  \label{clmcorrect} Assume $\mathsf{GCH}$.
  One can force $\mathfrak{gp}=\lambda$ whenever $\lambda>\cf(\lambda)=\omega_1$.
  Moreover, this statement has no consistency strength.
\end{claim}

\par\noindent\emph{Proof}. \newline
Let $\bar{\mathbb{S}}=\prod_{\alpha\in\lambda}\mathbb{S}_\alpha$ be the forcing notion of Abraham and Shelah from \cite{MR830084}.
Each component $\mathbb{S}_\alpha$ is a tow-step iteration in which the first step adds a Cohen real and the second step adds a club to $\aleph_1$.
For us it is important to observe that $\bar{\mathbb{S}}$ adds $\lambda$-many Cohen reals.
Let $G\subseteq\bar{\mathbb{S}}$ be $V$-generic.
In $V[G]$ one has a family $\mathcal{C}=\{C_\alpha:\alpha\in\lambda\}$ which witnesses the failure of the Galvin property.

Notice that $\cf(\lambda)>\omega_1$ is required in \cite{MR830084} only for the purpose of the ultimate failure, that is $2^{\omega_1}=\lambda$.
This point is elaborated in \cite{MR4611828}.
In our case, the family $\mathcal{C}$ is not a witness for the ultimate failure since $2^{\omega_1}>\lambda$ in $V[G]$ because $\cf(\lambda)=\omega_1$.
However, $2^\omega=\lambda$ in $V[G]$ since $\bar{\mathbb{S}}$ adds only $\lambda$-many Cohen reals and $\lambda^\omega=\lambda$ in the ground model.

Thus, we know that $\mathfrak{gp}\geq\lambda$ by the above forcing.
But we also have the cardinal arithmetic constellation of $2^\omega=\lambda<2^{\omega_1}$ and hence ${\rm Gal}(\mathscr{D}_{\aleph_1},\aleph_1,\lambda^+)$ holds as proved in \cite{MR3604115}.
Therefore, $\mathfrak{gp}=\lambda$ as required.

\hfill \qedref{clmcorrect}

One can force with \cite{MR830084} where $\lambda>\cf(\lambda)=\omega$ as well.
The same argument will show that $\mathfrak{gp}\geq\lambda$.
However, $2^\omega>\lambda$ in this case, so one cannot prove that $\mathfrak{gp}=\lambda$.
In fact, it seems that one can define explicitly a family of size $\lambda^+$, witnessing the failure of the Galvin property in this model.

We conclude this discussion with an interesting question raised by Eran Alouf, \cite{ea}.
In the above model for $\cf(\mathfrak{gp})=\omega_1$ there are a cardinal $\mu>\cf(\mu)=\omega_1$ and a family $\mathcal{D}\subseteq\mathscr{D}_{\aleph_1}$ of size $\mu$, such that $\bigcap{B}\notin\mathscr{D}_{\aleph_1}$ whenever $B\in[\mathcal{D}]^{\aleph_1}$.
This family is one step towards proving that $\mathfrak{gp}=\mu$.
There is, however, another scenario in which $\mathfrak{gp}=\mu$ and hence $\cf(\mathfrak{gp})=\omega_1$.
Namely, for every $\kappa<\mu$ there is a family $\mathcal{D}_\kappa$ witnessing the failure of the Galvin property, but there is no such of family of size $\mu$.

\begin{question}
  \label{qalouf} Is it consistent that $\mu>\cf(\mu)=\omega_1, \mathfrak{gp}=\mu$ and yet every $\mathcal{D}\in[\mathscr{D}_{\aleph_1}]^\mu$ satisfies the Galvin property?
\end{question}

Our last observation is related to the possibility of countable cofinality for the Galvin number.

\begin{remark}
  \label{rcountablecof} Suppose that $\mu>\cf(\mu)=\omega$ and $\mathfrak{gp}=\mu$.
  Then $2^\theta=\mu^+$ for every $\aleph_0\leq\theta\leq\mu$.
\end{remark}

\par\noindent\emph{Proof}. \newline
By a result of Shelah, if $2^\mu>\mu^+$ then there is a good scale at $\mu^+$, see \cite{MR1318912} or \cite[Section 19]{MR2160657}.
Thus, if $\mathfrak{gp}=\mu$ then necessarily $2^\mu=\mu^+$.
Likewise, if $2^{\aleph_0}<\mu$ then $\mathfrak{gp}\leq 2^{\aleph_0}<\mu$.
Hence, if $\mathfrak{gp}=\mu$ then $2^{\aleph_0}\geq\mu$ and then necessarily $2^{\aleph_0}\geq\mu^+$.
Since $2^{\aleph_0}\leq 2^\mu=\mu^+$, we are done.

\hfill \qedref{rcountablecof}

\newpage

\bibliographystyle{alpha}
\bibliography{arlist}

\begin{thebibliography}{GHHM22}

\bibitem[Alo24]{ea}
Eran Alouf.
\newblock Private communication.
\newblock 2024.

\bibitem[AM10]{MR2768693}
Uri Abraham and Menachem Magidor.
\newblock Cardinal arithmetic.
\newblock In {\em Handbook of set theory. {V}ols. 1, 2, 3}, pages 1149--1227.
  Springer, Dordrecht, 2010.

\bibitem[AS86]{MR830084}
U.~Abraham and S.~Shelah.
\newblock On the intersection of closed unbounded sets.
\newblock {\em J. Symbolic Logic}, 51(1):180--189, 1986.

\bibitem[BGP22]{bgp}
Tom Benhamou, Shimon Garti, and Alejandro Poveda.
\newblock Galvin's property at large cardinals and the axiom of choice.
\newblock {\em Israel Journal of Mathemtics}, (accepted), 2022.

\bibitem[BGP23]{MR4611828}
Tom Benhamou, Shimon Garti, and Alejandro Poveda.
\newblock Negating the {G}alvin property.
\newblock {\em J. Lond. Math. Soc. (2)}, 108(1):190--237, 2023.

\bibitem[BHM75]{MR0369081}
J.~E. Baumgartner, A.~H\c{a}j\c{n}al, and A.~Mate.
\newblock Weak saturation properties of ideals.
\newblock In {\em Infinite and finite sets ({C}olloq., {K}eszthely, 1973;
  dedicated to {P}. {E}rd\H{o}s on his 60th birthday), {V}ol. {I}}, pages
  137--158. Colloq. Math. Soc. J\'{a}nos Bolyai, Vol. 10. 1975.

\bibitem[Cum05]{MR2160657}
James Cummings.
\newblock Notes on singular cardinal combinatorics.
\newblock {\em Notre Dame J. Formal Logic}, 46(3):251--282, 2005.

\bibitem[Dic59]{tale}
Charles Dickens.
\newblock {\em A Tale of Two Cities}.
\newblock All The Year Round, A Weekly Journal. 1859.

\bibitem[EHR65]{MR0202613}
P.~Erd{\H{o}}s, A.~Hajnal, and R.~Rado.
\newblock Partition relations for cardinal numbers.
\newblock {\em Acta Math. Acad. Sci. Hungar.}, 16:93--196, 1965.

\bibitem[Gar17]{MR3604115}
Shimon Garti.
\newblock Weak diamond and {G}alvin's property.
\newblock {\em Period. Math. Hungar.}, 74(1):128--136, 2017.

\bibitem[Gar18]{MR3787522}
Shimon Garti.
\newblock Tiltan.
\newblock {\em C. R. Math. Acad. Sci. Paris}, 356(4):351--359, 2018.

\bibitem[Gar20]{MR4101445}
Shimon Garti.
\newblock Polarized relations at singulars over successors.
\newblock {\em Discrete Math.}, 343(9):111961, 9, 2020.

\bibitem[GHHM22]{MR4423479}
Shimon Garti, Yair Hayut, Haim Horowitz, and Menachem Magidor.
\newblock Forcing axioms and the {G}alvin number.
\newblock {\em Period. Math. Hungar.}, 84(2):250--258, 2022.

\bibitem[Hen83]{MR722169}
J.~M. Henle.
\newblock Magidor-like and {R}adin-like forcing.
\newblock {\em Ann. Pure Appl. Logic}, 25(1):59--72, 1983.

\bibitem[Kle77]{MR0479903}
Eugene~M. Kleinberg.
\newblock {\em Infinitary combinatorics and the axiom of determinateness}.
\newblock Lecture Notes in Mathematics, Vol. 612. Springer-Verlag, Berlin-New
  York, 1977.

\bibitem[LS67]{MR224458}
A.~L\'{e}vy and R.~M. Solovay.
\newblock Measurable cardinals and the continuum hypothesis.
\newblock {\em Israel J. Math.}, 5:234--248, 1967.

\bibitem[Mag78]{MR465868}
Menachem Magidor.
\newblock Changing cofinality of cardinals.
\newblock {\em Fund. Math.}, 99(1):61--71, 1978.

\bibitem[She94]{MR1318912}
Saharon Shelah.
\newblock {\em Cardinal arithmetic}, volume~29 of {\em Oxford Logic Guides}.
\newblock The Clarendon Press, Oxford University Press, New York, 1994.
\newblock Oxford Science Publications.

\bibitem[Ste10]{MR2768698}
John~R. Steel.
\newblock An outline of inner model theory.
\newblock In {\em Handbook of set theory. {V}ols. 1, 2, 3}, pages 1595--1684.
  Springer, Dordrecht, 2010.

\bibitem[Wil77]{MR3075383}
Neil~H. Williams.
\newblock {\em Combinatorial set theory}, volume~91 of {\em Studies in Logic
  and the Foundations of Mathematics}.
\newblock North-Holland Publishing Co., Amsterdam, 1977.

\end{thebibliography}

\end{document}